\documentclass[12pt]{amsart}
\usepackage {amssymb,amsmath}
\makeatletter
\def\@cite#1#2{{\m@th\upshape\bfseries%
[{#1\if@tempswa{\m@th\upshape\mdseries, #2}\fi}]}} \makeatother
\newtheorem{thm}{Theorem}[section]
\newtheorem{lem}[thm]{Lemma}

\newtheorem{defn}[thm]{Definition}

\newcommand{\Prf}{\noindent\textbf{Proof.\ }}
\newcommand{\bx}{\strut\hfill$\blacksquare$\medbreak}

\newcommand{\ca}{\mathrm{C}^*}

\newenvironment{spmatrix}{\left(\begin{smallmatrix}}{\end{smallmatrix}\right)}

\newcommand{\bbC}{{\mathbb{C}}}

 \newcommand{\B}{{\mathcal{B}}}
 \newcommand{\C}{{\mathcal{C}}}
 
 \newcommand{\E}{{\mathcal{E}}}

\renewcommand{\H}{{\mathcal{H}}}

 \newcommand{\K}{{\mathcal{K}}}

\renewcommand{\P}{{\mathcal{P}}}
 
 \newcommand{\R}{{\mathcal{R}}}
\newcommand{\upchi}{{\raise.35ex\hbox{$\chi$}}}
\newcommand{\fA}{{\mathfrak{A}}}

\newcommand{\qforal}{\quad\text{for all}\quad}

\newcommand{\spn}{\operatorname{span}}

\newcommand{\Tr}{\operatorname{Tr}}

\def\bra#1{\langle #1|}
\def\ket#1{|#1 \rangle}
\def\kb#1#2{|#1\rangle\!\langle #2 |}
\def\one{{\mathchoice{\rm 1\mskip-4mu l}{\rm 1\mskip-4mu l}{\rm 1\mskip-4.5mu l}{\rm
1\mskip-5mu l}}}

\begin{document}

%%%%%%%%%%%%%%%%%%%%%%%%%%%%%%%%%%%%%%%%%%%%%%
%%%%%%%%%%
\title[Operator Quantum Error Correction]{A Brief Introduction To Operator Quantum Error
Correction\\  \vspace{0.5cm} \tiny {\it Dedicated to John Holbrook
on the occasion of his 65th birthday}}
\author[D.W. Kribs]{David W. Kribs}
\thanks{2000 {\it Mathematics Subject Classification.} 47L90,  81P68.}
\thanks{{\it key words and phrases.} completely positive map, quantum operation, error model, generalized noiseless subsystem,
decoherence-free subspace, operator quantum error correction.}

\address{Department of Mathematics and Statistics,
University of Guelph, Guelph, Ontario, Canada  N1G 2W1}
\address{Institute for Quantum Computing, University of
Waterloo, Waterloo, ON, CANADA N2L 3G1}

%\date{}
%
\begin{abstract}
We give a short introduction to operator quantum error correction.
This is a new protocol for error correction in quantum computing
that has brought the fundamental methods under a single umbrella,
and has opened up new possibilities for protecting quantum
information against undesirable noise. We describe the various
conditions that characterize correction in this scheme.
\end{abstract}
\maketitle
%%%%%%%%%%%%%%%%%%%%%%%%%%%%%%%%%%%%%%%%%%%%%%
%%%%%%%%%%%

%%%%%%%%%%%%%%%%%%%%%%%%%%%%%%%%%%%%%%
\section{Introduction}\label{S:intro}
%%%%%%%%%%%%%%%%%%%%%%%%%%%%%%%%%%%%%%

In this paper we give an introduction to some of the mathematical
aspects of quantum error correction, with an emphasis on the
unified approach -- called {\it operator quantum error correction}
-- recently introduced in \cite{KLP05,KLPL05}. The field of
quantum error correction took flight during the mid 1990's
\cite{Sho95a,Ste96a,BDSW96a,KL97a}. A central goal of this young
field is to help construct quantum computers via the development
of schemes that allow for the protection of quantum information
against the noise associated with evolution of quantum systems. As
it turns out, many of the problems in quantum error correction
have an operator theoretic flavour. Here we shall briefly discuss
the fundamental error correction protocols in quantum computing.
We also describe various conditions that characterize correction
in Operator QEC, and provide a new operator proof for the main
testable condition in this scheme.

First let us briefly discuss the basic setting for quantum
computation. See \cite{NC00,Kri05a} as examples of more extensive
introductions. To each quantum system, the postulates of quantum
mechanics associate a Hilbert space $\H = \bbC^n$. The finite
dimensional case is the current focus in quantum computing for
experimental reasons. A two-level quantum system is represented on
$\H = \bbC^2$. This could describe, for instance, the ground and
excited energy states of an electron in an atom. These are the two
classical states that we observe, corresponding to an orthonormal
basis $\{\ket{0},\ket{1}\}$ for $\bbC^2$. However, quantum
mechanics dictates that any linear combination $\ket{\psi}$ of
these classical states is an allowable state, even though we only
observe either the $\ket{0}$ or $\ket{1}$ state. When it is a
non-trivial linear combination, $\ket{\psi}$ is said to be in a
{\it superposition} of the classical states. A unit vector
$\ket{\psi}$ is the fundamental quantum bit of information, also
called a ``qubit''. Equivalently, we could consider the rank one
projection $\kb{\psi}{\psi}$. The corresponding $n$-qubit
composite system is realized on $\H = (\bbC^2)^{\otimes n} =
\bbC^{2^n}$ with orthonormal basis $\big\{\ket{i_1\cdots
i_n}\equiv \ket{i_1}\otimes\ldots\otimes\ket{i_n}:
i_j\in\{0,1\}\big\}$ determined by the underlying two-level
system.

More generally, we will only know that our system is in one of
several states with various possibilities. So the direct
generalization of a classical probability distribution in quantum
information theory is a positive matrix $\rho$ with trace equal to
one, a so-called {\it density matrix}.

A fundamental problem in quantum computation is to physically
manipulate the superpositions inherent in quantum systems, without
collapsing or ``decohering'' them. To accomplish this, methods
must be developed to correct the errors that occur as quantum
information is transferred from one physical location to the next
inside, for instance, a quantum computer. To deal with this
problem we must discuss evolution of quantum systems, a subject to
which we now turn.

%%%%%%%%%%%%%%%%%%%%%%%%%%%%%%%%%%%%%%
\section{Evolution of Quantum Systems and Error Correction}\label{S:evolve}
%%%%%%%%%%%%%%%%%%%%%%%%%%%%%%%%%%%%%%

The reversibility postulate of quantum mechanics implies that
evolution in a closed quantum system occurs via unitary maps. From
the discrete perspective, if we take a snapshot of this evolution,
then a density matrix $\rho$ will encode the possible states of
the system with various probabilities at a given time. An
evolution of the system corresponds to a map $\rho \mapsto U\rho
U^\dagger$ for some unitary operator $U$.

In the context of quantum computing, the quantum systems of
interest are ``open'' as they are exposed to external environments
during computations. In such cases, the open system is regarded as
part of a larger closed quantum system given by the composite of
the system and the environment. If $\H_S$ and $\H_E$ are the
system and environment Hilbert spaces, then the closed system is
represented on $\H = \H_E \otimes \H_S$.

The characterization of evolution in open quantum systems requires
first that density operators are mapped to density operators; i.e.
probability densities are mapped to probability densities. Thus,
such a map must be positive and trace preserving. However, this
property must be preserved when the system is exposed to all
possible environments. In terms of the map, if $\E$ describes an
evolution of the system, then the map $id_E\otimes
\E:\B(\H_E\otimes\H_S)\rightarrow\B(\H_E\otimes\H_S)$ must also be
positive and trace preserving for all $E$. Hence, the widely
accepted working definition of a {\it quantum operation} (or
evolution, or channel) on a Hilbert space $\H$, is a completely
positive, trace preserving map  $\E$ on $\B(\H)$ (CPTP for short).

Deriving from a theorem of Choi \cite{Cho75} and Kraus
\cite{Kra71}, every CPTP map $\E: \B(\H)\rightarrow\B(\H)$ has an
``operator-sum representation'' of the form $\E(\rho) = \sum_a E_a
\rho E_a^\dagger$ for some set of (non-unique) operators
$\{E_a\}\subseteq\B(\H)$ with $\sum_a E_a^\dagger E_a = \one$. The
$E_a$ are called the {\it noise operators} or {\it errors}
associated with $\E$. In the context of quantum error correction,
it is precisely the effects of these errors that must be
mitigated. As a short hand, we write $\E = \{E_a\}$ when an error
model for $\E$ is known.

Error correction in quantum computing is a much more delicate
problem in comparison  to its classical counterpart. As a simple
observation, consider that the only errors that occur classically
are some version of bit flips; e.g., $\ket{0}$ goes to $\ket{1}$
or vice-versa. More generally, in quantum computing subtleties
arise from the fact that a given qubit can be corrupted to an {\it
infinite} number of possible superpositions. In terms of operators
on single bits or qubits for instance, whereas the Pauli bit flip
matrix $X =
\begin{spmatrix} 0 & 1 \\ 1 & 0 \end{spmatrix}$ is the fundamental
classical error matrix, {\it any} unitary matrix is a possible
error in quantum computing. Of course, there are many other
issues, such as the fabled ``No Cloning Theorem''. The linearity
of quantum mechanics implies that the analogue of the classically
well-used repetition code does not extend to arbitrary qubits
$\ket{\psi}\mapsto\ket{\psi}\otimes\ket{\psi}$. Fortunately,
methods have been, and are being, developed to overcome these
challenges.

%%%%%%%%%%%%%%%%%%%%%%%%%%%%%%%%%%%%%%
\section{Standard Model of Quantum Error
Correction}\label{S:standardmodel}
%%%%%%%%%%%%%%%%%%%%%%%%%%%%%%%%%%%%%%

The ``Standard Model'' of quantum error correction
\cite{Sho95a,Ste96a,BDSW96a,KL97a} involves triples $(\R,\E,\C)$
where $\C$ is a subspace, a {\it quantum code}, of a Hilbert space
$\H$ associated with a given quantum system, and the {\it error}
$\E$ and {\it recovery} $\R$ are quantum operations on $\B(\H)$.

Recall from the discussions above that we are forced by quantum
mechanics to consider subspaces $\C$ as sets of codes, as linear
combinations of classical codewords are perfectly allowable
codewords in this setting. In the trivial case, when $\E=\{U\}$ is
implemented by a single unitary error operator, the recovery is
just the reversal operation $\R = \{U^\dagger\}$; that is,
\[
\rho \,\,\stackrel{\E}{\longrightarrow} \,\, U \rho U^\dagger
\stackrel{\R}{\longrightarrow}\,\, U^\dagger (U\rho U^\dagger) U =
\rho.
\]
Of course, here there is no need to restrict the input operators
$\rho$.

More generally, the set $(\R,\E,\C)$ forms an ``error triple'' if
$\R$ undoes the effects of $\E$ on $\C$ in the following sense:
\begin{eqnarray}\label{reverse}
(\R\circ \E) \, (\sigma) = \sigma \quad \forall\, \sigma
\in\B(\C),
\end{eqnarray}
where $\C$ is naturally regarded as embedded inside $\H$.

When there exists such  an $\R$ for a given pair $\E,\C$, the
subspace $\C$ is said to be {\it correctable for $\E$}. The
existence of a recovery operation $\R$ of $\E=\{E_a\}$ on $\C$ is
characterized by the following condition \cite{BDSW96a,KL97a}:
There exists a scalar matrix $\Lambda = (\lambda_{ab})$ such that
\begin{equation}
P_\C E_a^\dagger E_b P_\C = \lambda_{ab}P_\C \quad \forall\, a,b,
\label{eq:standard}
\end{equation}
where $P_\C$ is the projection of $\H$ onto $\C$. It is not hard
to see that this condition is independent of the operator-sum
representation for $\E$. We note that Eq.~(\ref{eq:standard}) is a
special case of Eq.~(\ref{condition}) below.

The motivating case of an error model that satisfies
Eq.~(\ref{eq:standard}) occurs when the restrictions
$E_a|_{P_\C\H}= E_a|_\C$ of the noise operators to $\C$ are scalar
multiples of unitary operators $U_a$, such that the subspaces
$U_a\C$ are mutually orthogonal. In this situation the positive
scalar matrix $\Lambda$ is diagonal. A correction operation here
may be constructed by an application of the measurement operation
determined by the subspaces $U_a\C$, followed by the reversals of
the corresponding restricted unitaries $U_a P_\C$. Specifically,
if $P_a$ is the projection of $\H$ onto $U_a\C$, then $\R =
\{U_a^\dagger P_a\}$ satisfies Eq.~(\ref{reverse}) for $\E$ on
$\C$.

Let us discuss a simple example. Let $\C$ be the subspace of $\H =
\bbC^8$ given by $\C = \spn\{\ket{000},\ket{111}\}$. Let $\E = \{
\frac{1}{\sqrt{3}} X_k : k = 1,2,3\}$ with the Pauli matrix $X$
and $X_1 = X\otimes \one_2 \otimes \one_2$, and similarly for
$X_2$, $X_3$. In this case, $\Lambda = \frac{1}{3} \one_3$. The
correction operation $\R$ may be constructed as above.

%%%%%%%%%%%%%%%%%%%%%%%%%%%%%%%%%%%%%%%%%%%%%%%%%%%%%%%%%%%%%
\section{Noiseless Subsystems}\label{S:genNS}
%%%%%%%%%%%%%%%%%%%%%%%%%%%%%%%%%%%%%%%%%%%%%%%%%%%%%%%%%%%%%

To describe the notion of noiseless subsystems from
\cite{KLP05,KLPL05}, we begin with a decomposition of the system
Hilbert space
\begin{equation*}
\H = \bigoplus_J \H^A_J\otimes\H^B_J,
\end{equation*}
where the ``noisy subsystems" $\H^A_J$ have dimension $m_J$ and
the ``noiseless subsystems" $\H^B_J$ have dimension $n_J$. We
focus on the case where information is encoded in a single
noiseless sector of $\B(\H)$, so
\begin{equation*}
\H = (\H^A \otimes \H^B) \oplus \K \label{eq:decomp}
\end{equation*}
with $\dim(\H^A) = m$, $\dim(\H^B) = n$ and $\dim \K=\dim\H - mn$.
We shall write $\sigma^A$ for operators in $\B(\H^A)$ and
$\sigma^B$ for operators in $\B(\H^B)$.

Let $\{\ket{\alpha_k}:1\leq k \leq m\}$ be an orthonormal basis
for $\H^A$ and let
\begin{equation*}
\{ P_{kl} = \kb{\alpha_k}{\alpha_l}\otimes\one_{n}: 1\leq k,l \leq
m\}
\end{equation*}
be the corresponding family of matrix units in $\B(\H^A)\otimes
\one^B$. Recall that the partial trace over $A$ on
$\H^A\otimes\H^B$ is the quantum operation defined on elementary
tensors by $\Tr_A(\sigma^A\otimes\sigma^B) =
\Tr(\sigma^A)\sigma^B$.

Define for a fixed decomposition $\H = (\H^A\otimes\H^B) \oplus\K$
the operator semigroup
\begin{eqnarray}\label{eq:semigroup}
\fA = \{\sigma\in\B(\H) : \sigma = \sigma^A\otimes\sigma^B,\,{\rm
for\,\,some}\,\, \sigma^A\,{\rm and}\, \sigma^B\}.
\end{eqnarray}
For notational purposes, we assume that bases have been chosen and
define the matrix units $P_{kl}$ as above, so that $P_k = P_{kk}$,
$P_\fA \equiv P_1+\ldots + P_m$ and $P_\fA\H = \H^A\otimes\H^B$.
We also define a map $\P_\fA$ by the action $\P_\fA(\cdot) =
P_\fA(\cdot) P_\fA$. The following result motivates the
(generalized) definition of NS's from \cite{KLP05,KLPL05}. (See
\cite{KLPL05} for a proof.)

\begin{lem}
Given a fixed decomposition $\H = (\H^A\otimes\H^B) \oplus\K$ and
a quantum operation $\E$ on $\B(\H)$, the following three
conditions are equivalent:
\begin{enumerate}
\item[({\it 1})] $\forall\sigma^A\ \forall\sigma^B,\ \exists
\tau^A\ :\ \E(\sigma^A\otimes\sigma^B) = \tau^A\otimes\sigma^B$
\item[({\it 2})] $ \forall\sigma^B,\ \exists \tau^A\ :\
\E(\one^A\otimes\sigma^B) = \tau^A \otimes \sigma^B$ \item[({\it
3})] $\forall\sigma\in \fA\ :\ \big(\Tr_A\circ \P_\fA\circ
\E\big)(\sigma) =\Tr_A(\sigma)$.
\end{enumerate}
\label{lemma:generalNS}
\end{lem}

\begin{defn}\label{defn:NS}
{\rm The $\H^B$ sector of the semigroup $\fA$ encodes a  {\it
noiseless subsystem} for $\E$ when it satisfies the equivalent
conditions of Lemma~\ref{lemma:generalNS}.}
\end{defn}

The NS framework discussed here is a generalization of both the
``Standard NS'' \cite{KLV00a,Zan01b,KBLW01a} and
``Decoherence-Free Subspace'' \cite{PSE96,DG97c,ZR97c,LCW98a}
methods of passive error correction, both of which are used for
unital quantum operations. The method described here applies to
{\it all} CPTP maps. See \cite{KLP05,KLPL05} for more discussions
on this point.

As a simple example of how such subsystems naturally arise, let
$\Phi:\B(\H^A)\rightarrow\B(\H^A)$ be an arbitrary CPTP map and
let $\Psi:\B(\H^B)\rightarrow\B(\H^B)$ be CPTP with a Standard NS
$\H_0^B\subseteq\H^B$; i.e., $\Psi(\rho) = \rho$ for all
$\rho\in\B(\H_0^B)$. Then $\H_0^B$ encodes a noiseless subsystem
inside $\H^A\otimes\H^B$ for the map $\E = \Phi\otimes\Psi:
\B(\H^A\otimes\H^B)\rightarrow\B(\H^A\otimes\H^B)$.

To be of use in practical applications, we need testable
conditions for a map $\E = \{E_a\}$ to admit a NS described by a
semigroup $\fA$. Towards this end, we have proved the following
theorem.

\begin{thm}\label{thm:NS}
Let $\E = \{E_a\}$ be a quantum operation on $\B(\H)$ and let
$\fA$ be a semigroup in $\B(\H)$ as above. Then the following
three conditions are equivalent:
\begin{itemize}
\item[({\it 1})] The $\H^B$ sector of $\fA$ encodes a noiseless
subsystem for $\E$. \item[({\it 2})] The subspace $P_\fA \H = \H^A
\otimes \H^B$ is invariant for the operators $E_a$ and the
restrictions $E_a|_{P_\fA\H}$ belong to the algebra
$\B(\H^A)\otimes\one^B$. \item[({\it 3})] The following two
conditions hold:
\begin{equation}
P_k E_a P_l = \lambda_{akl} P_{kl} \quad\forall\, a,k,l
\label{eq:cond1}
\end{equation}
for some set of scalars $(\lambda_{akl})$ and
\begin{equation}
 E_a P_\fA = P_\fA E_a P_\fA \quad\forall\, a. \label{eq:cond2}
\end{equation}
\end{itemize}
\end{thm}

\Prf Since the matrix units $\{P_{kl}\}$ generate $\B(\H^A)\otimes
\one^B$, it follows that {\it (3)} is a restatement of {\it (2)}.
Here we sketch the proof of the equivalence of {\it (1)} and {\it
(3)}, see \cite{KLP05,KLPL05} for details. To prove the necessity
of Eqs.~(\ref{eq:cond1}), (\ref{eq:cond2}) for {\it (1)}, it
follows from properties of the map $\Gamma = \{P_{kl}\}$ and
Lemma~\ref{lemma:generalNS} that there exist scalars
$\mu_{kiajl,k'l'}$ such that
\begin{equation}
P_{ki}E_aP_{jl} = \sum_{k'l'} \mu_{kiajl,k'l'} P_{k'l'}.
\end{equation}
Multiplying both sides of this equality on the  right by $P_l$ and
on the left by $P_k$, we see that $\mu_{kiajl,k'l'} = 0$ when
$k\neq k'$ or $l\neq l'$. This implies Eq.~(\ref{eq:cond1}) with
$\lambda_{akl} = \mu_{kkall,kl}$. Equation~(\ref{eq:cond2})
follows from Lemma~\ref{lemma:generalNS} and consideration of the
operator-sum representation for $\E$.

On the other hand, if Eqs.~(\ref{eq:cond1}), (\ref{eq:cond2})
hold, then for all $\sigma= P_\fA \sigma \in\fA$ we have
\[
\E(\sigma) = \sum_{a,k,k'} P_kE_a\sigma E_a^\dagger P_{k'}.
\]
This implies that for all $\sigma =
\sigma^A\otimes\sigma^B\in\fA$,
\begin{eqnarray*}
\E(\sigma^A\otimes\sigma^B) &=& \sum_{a,k,k',l,l'} P_kE_aP_l(\sigma^A\otimes\sigma^B) P_{l'} E^\dagger_a P_{k'} \\
&=& \sum_{a,k,k',l,l'} \lambda_{akl}\overline{\lambda}_{ak'l'}
P_{kl} (\sigma^A\otimes\sigma^B) P_{l'k'}.
\end{eqnarray*}
Condition {\it (1)} now follows from the fact that the matrix
units $P_{kl}$ act trivially on the $\B(\H^B)$ sector.
 \bx

%%%%%%%%%%%%%%%%%%%%%%%%%%%%%%%%%%%%%%%%%%%%%%
%%%%%
\section{Operator Quantum Error Correction}\label{S:unified}
%%%%%%%%%%%%%%%%%%%%%%%%%%%%%%%%%%%%%%%%%%%%%%
%%%%%%%%%%%%%

The Operator QEC approach consists of  triples $(\R,\E,\fA)$ where
$\R$ and $\E$ are quantum operations on some $\B(\H)$, and $\fA$
is a semigroup in $\B(\H)$ defined as above with respect to a
fixed decomposition $\H = (\H^A \otimes \H^B) \oplus \K$.

\begin{defn}
{\rm Given a triple $(\R,\E,\fA)$ we say that the $\H^B$ sector of
$\fA$ is {\it correctable for $\E$} if
\begin{eqnarray}\label{newid}
\big(\Tr_A \circ \P_\fA \circ\R \circ \E \big) (\sigma) =
\Tr_A(\sigma) \qforal \sigma \in \fA.
\end{eqnarray}
}
\end{defn}

Equivalently, $(\R,\E,\fA)$ is a correctable triple if the $\H^B$
sector of the semigroup $\fA$ encodes a noiseless subsystem for
the error map $\R\circ\E$. Table~1 indicates the special cases
captured by Operator QEC. Our choice of terminology here is
motivated by the fact that correctable codes in this scheme take
the form of operator algebras and operator semigroups. We point
the reader to \cite{KLPL05,Bac05} for examples of error triples on
subsystems that require non-trivial recovery operations, and
\cite{NP05,LS05} for other recent related work.

\begin{table}

Table 1: Special Cases of Operator QEC \vspace{0.1in}

\begin{tabular}{||c||c||} \hline\hline
$\fA = \mbox{subspace}$ & Standard QEC \\ \hline $\R = id$ &
(Generalized) NS \\ \hline $\R = id \,\, + \,\, \fA =
\mbox{algebra}$ & Standard NS \\ \hline
$\R = id \,\, + \,\, \fA = \mbox{subspace}$ & DFS \\
\hline\hline
\end{tabular}

\end{table}

An important feature of Operator QEC is that a semigroup $\fA$ is
correctable exactly when the $\ca$-algebra $\fA_0 = \one^A \otimes
\B(\H^B)$ can be corrected precisely.

\begin{thm}\label{thm:opalgequiv}
Let $\E = \{E_a\}$ be a quantum operation on $\B(\H)$ and let
$\fA$ be a semigroup in $\B(\H)$ as above. Then the $\H^B$ sector
of $\fA$ is correctable for $\E$ if and only if there is a quantum
operation $\R$ on $\B(\H)$ such that
\begin{eqnarray}\label{opalgcorrect}
(\R\circ\E)(\sigma) = \sigma \quad \forall\, \sigma \in \one^A
\otimes \B(\H^B).
\end{eqnarray}
\end{thm}

\Prf If Eq.~(\ref{opalgcorrect}) holds, then condition~{\it (2)}
of Lemma~\ref{lemma:generalNS} holds for $\R\circ\E$ with $\tau^A
= \one^A$ and hence $\fA$ is correctable for $\E$. On the other
hand, suppose that $\fA$ is correctable for $\E$ and
condition~{\it (2)} of Lemma~\ref{lemma:generalNS} holds for
$\R\circ\E$. Note that the map $\Gamma' =
\{\frac{1}{\sqrt{m}}P_{kl}\}$ is trace preserving on $\B(\H^A
\otimes \H^B)$. Thus, from basic properties of the map $\Gamma =
\{P_{kl}\}$, we have for all $\sigma^B$,
\begin{eqnarray}\label{opalgequiveqn}
(\Gamma' \circ\R\circ\E) (\one^A \otimes \sigma^B) = \Gamma'
(\tau^A \otimes \sigma^B) \propto \one^A \otimes \sigma^B.
\end{eqnarray}
By trace preservation the proportionality factor must be one, and
hence Eq.~(\ref{opalgcorrect}) is satisfied for $(\Gamma'
\circ\R)\circ\E$. The map $\Gamma^\prime$ may be extended to a
quantum operation on $\B(\H)$ by including the projection
$P_\fA^\perp$ onto $\K$ as a noise operator. As this does not
effect the calculation (\ref{opalgequiveqn}), the result follows.
 \bx

We next give a testable condition, Eq.~(\ref{condition}), that
characterizes correction in the Operator QEC regime. Notice that
this is a generalization of Eq.~({\ref{eq:standard}) for Standard
QEC. This condition was introduced in \cite{KLP05} and necessity
was established. Sufficiency was proved in \cite{KLPL05} up to a
set of technical conditions, and more recently in \cite{NP05} with
full generality. (The work of \cite{NP05} also links this
condition with an interesting information theoretic condition.)
Here we include a sketch of the proof of necessity from
\cite{KLPL05}, and a new operator theoretic version of the proof
of sufficiency sketched in \cite{NP05}. We assume that matrix
units $\{P_{kl}= \kb{\alpha_k}{\alpha_l} \otimes \one^B \}$ inside
$\B(\H^A)\otimes\one^B$ have been chosen as above.

\begin{thm}\label{unifiedthm}
Let $\E = \{E_a\}$ be a quantum operation on $\B(\H)$ and let
$\fA$ be a semigroup in $\B(\H)$ as above. For the $\H^B$ sector
of $\fA$ to be correctable for $\E$, it is necessary and
sufficient that there are scalars $\Lambda = ( \lambda_{abkl} )$
such that
\begin{eqnarray}\label{condition}
P_{k} E_a^\dagger E_b P_l = \lambda_{abkl} P_{kl} \quad \forall\,
a,b,k,l.
\end{eqnarray}
\end{thm}

\Prf For necessity, note first that Theorem~\ref{thm:opalgequiv}
gives us a CPTP map $\R$ on $\B(\H)$ such that $\R\circ\E$ acts as
the identity channel on $\fA_0 = \one^A \otimes \B(\H^B)\subseteq
\B(\H)$.

Suppose that  $\R = \{ R_b \}$. The noise operators for the
operation $\R \circ \E$ are $\{R_b E_a\}$, and using arguments
similar to those of Theorem~\ref{thm:NS} (see \cite{KLPL05} for
details) we may find scalars $\mu_{abkl}$ such that
\[
P_k R_b E_a P_l = \mu_{abkl} P_{kl} \quad \forall a,b,k,l.
\]
Consider the products
\begin{eqnarray*}
\big(P_k R_b E_a P_l \big)^\dagger \big(P_{k'} R_b E_{a'}
P_{l'}\big) &=& \big( \overline{\mu_{abkl}} P_{lk} \big) \big(
\mu_{a'bk'l'} P_{k'l'} \big) \\
&=& \left\{ \begin{array}{cl} (\overline{\mu_{abkl}}\mu_{a'bkl'}
)P_{ll'} & \mbox{if $k=k'$} \\
0 & \mbox{if $k\neq k'$}
\end{array}\right. .
\end{eqnarray*}
Now, the subspace $\C$ can be shown to be invariant for the noise
operators $R_b E_a$. Hence for fixed $a, a'$ and $l,l'$ we use
$\sum_b R_b^\dagger R_b = \one$ to obtain
\begin{eqnarray*}
\Big( \sum_{b,k} \overline{\mu_{abkl}} \mu_{a'bkl'} \Big) P_{ll'}
&=& \sum_{b,k} \big(P_l E_a^\dagger R_b^\dagger P_k \big) \big(P_k
R_b E_{a'} P_{l'} \big) \\
&=& \sum_b P_l E_a^\dagger R_b^\dagger P_\fA R_b E_{a'}
P_{l'} \\
&=&  P_l E_a^\dagger \Big( \sum_b R_b^\dagger R_b \Big) E_{a'} P_{l'} \\
&=& P_l E_a^\dagger E_{a'} P_{l'}
\end{eqnarray*}
The proof of necessity is completed by setting $\lambda_{aa'll'}
=\sum_{b,k} \overline{\mu_{abkl}} \mu_{a'bkl'}$ for all $a,a'$ and
$l,l'$.

For sufficiency, let us assume that Eq.~(\ref{condition}) holds.
Let $\sigma_k =\kb{\alpha_k}{\alpha_k}\in\B(\H^A)$, for $1\leq k
\leq m$, and define a CPTP map $\E_k:\B(\H^B)\rightarrow\B(\H)$ by
$ \E_k(\rho^B) \equiv \E(\sigma_k\otimes\rho^B). $ With $P\equiv
P_\fA$ and $E_{a,k} \equiv E_a P \ket{\alpha_k}$, it follows that
$\E_k = \{ E_{a,k}\} $. We shall find a CPTP map that globally
corrects {\it all} of the errors  $E_{a,k}$.

To do this, first note that  we may define a CPTP map
$\E_B:\B(\H^B)\rightarrow\B(\H)$ with error model
\[
\E_B = \big\{ \frac{1}{\sqrt{m}} E_{a,k} : \forall a, \, \forall
1\leq k \leq m\big\}.
\]
Then Eq.~(\ref{condition}) and $P=\sum_k P_k$ give us
\begin{eqnarray*}
\one^B E_{a,k}^\dagger E_{b,l} \one^B &=& \one^B \bra{\alpha_k} P
E_a^\dagger E_b P \ket{\alpha_l} \one^B \\
&=& \sum_{k',l'} \one^B \bra{\alpha_k} P_{k'} E_a^\dagger E_b
P_{l'} \ket{\alpha_l} \one^B \\
&=& \sum_{k',l'} \lambda_{abk'l'} \, \one^B \bra{\alpha_k}
P_{k'l'} \ket{\alpha_l} \one^B \\
&=&  \lambda_{abkl}  \one^B.
\end{eqnarray*}
In particular, Standard QEC implies the existence of a CPTP map
$\R: \B(\H)\rightarrow\B(\H^B)$ such that $(\R\circ\E_B)(\rho^B) =
\rho^B$ for all $\rho^B$.

This implies that
\begin{eqnarray*}
(\R\circ\E)(\one^A\otimes\rho^B) &=& \R \Big( \sum_k \E_k(\rho^B)
\Big) \\ &=& m \, \R \Big( \sum_{k,a} \frac{1}{m} E_{a,k} \rho^B
E_{a,k}^\dagger \Big) \\ &=& m \, \R \circ \E_B (\rho^B) = m
\rho^B.
\end{eqnarray*}
Hence we may define a CPTP ampliation map $I_\fA :
\B(\H^B)\rightarrow \B(\H)$ via $I_\fA(\rho^B) = \frac{1}{m}
(\one^A \otimes \rho^B$). Thus on defining $\R' \equiv I_\fA \circ
\R$, we obtain
\[
\big( \R' \circ\E\big)(\one^A\otimes\rho^B) = \one^A\otimes\rho^B
\quad \forall\, \rho^B\in\B(\H^B).
\]
The result now follows from an application of
Theorem~\ref{thm:opalgequiv}.
 \bx

%\begin{rem}
%{\rm Notice that (\ref{condition}) may be stated in matrix form as
% (fill in)
% }
%\end{rem}

%There is an interesting operator algebra consequence of the proof.
%[CHECK]
%
%\begin{cor}
%Let $\R = \{R_b\}$ be an error correction operation for $\E=
%\{E_a\}$ and $\fA$ as above. Then $\P_\fA\H$ is invariant for each
%$R_b E_a$  and the operators $R_b E_a P_\fA = P_\fA R_b E_a P_\fA$
%belong to the algebra  $\one^A \otimes \B(\H^B)$ for all $a,b$.
%\end{cor}

%%%%%%%%%%%%%%%%%%%%%%%%%%%%%%%%%%%%%%%%%%%%%%%%%%%%%%%%%%%
\section{Concluding Remark}\label{S:conclusion}
%%%%%%%%%%%%%%%%%%%%%%%%%%%%%%%%%%%%%%%%%%%%%%%%%%%%%%%%%%%%

The focus of research in quantum error correction has mainly been
on finite dimensional problems to this point. Primarily this
reflects the current status of experimental efforts to build
quantum computers, and the fact that many scientists working in
the area are closely linked with experimentalists. Thus, in the
author's opinion, there is an opportunity here for operator
theorists. In particular, mathematicians working in the field
have, for the most part, not had the luxury of exploring infinite
dimensional aspects and extensions of the quantum error correction
framework. It is expected that problems of this nature will
eventually be of experimental interest, and we expect they would
be of current mathematical interest.

%%%%%%%%%%%%%%%%%%%%%%%%%%%%%%%%%%%%%%%%%%%%%%
%%%%%%%%%%%%%
%%%%%%%%%%%%%%%%%%%%%%%%%%%%%%%%%%%%%%%%%%%%%%
%%%%%%%%%%%%%

\vspace{0.1in}

{\noindent}{\it Acknowledgements.} This paper was prepared as part
of the Proceedings of the 25th Anniversary Meeting of the Great
Plains Operator Theory Symposium held at the University of Central
Florida in June 2005. We thank the organizers and participants for
a  stimulating meeting. Thanks also to Palle Jorgensen for helpful
comments on an early draft of this paper.  This work was partially
supported by NSERC.

%%%%%%%%%%%%%%%%%%%%%%% REFERENCES
%%%%%%%%%%%%%%%%%%%%%%%%%%%%

%\newpage

%\begin{tabbing}
%{\it E-mail address}:xx\= \kill
%\noindent {\footnotesize\it Addresses}:
%\>{\footnotesize\sc Department of Mathematics}\\
%\>{\footnotesize\sc University of Iowa}\\
%\>{\footnotesize\sc Iowa City, IA\quad 52242}\\
%\>{\footnotesize\sc USA}\\
%\\
%\>{\footnotesize\sc Department of Mathematics and Statistics}\\
%\>{\footnotesize\sc Lancaster University}\\
%\>{\footnotesize\sc Lancaster, England}\\
%\>{\footnotesize\sc LA1 4YW}\\
%\\
%{\footnotesize\it E-mail addresses}:
%\>{\footnotesize\sf dkribs@math.uiowa.edu}\\
%\>{\footnotesize\sf s.power@lancaster.ac.uk}
%
%\end{tabbing}

\end{document}